\documentclass[10pt,letterpaper]{amsart}

\usepackage{mathrsfs}
\usepackage{latexsym}
\usepackage{amsmath,amsfonts,amsthm,amssymb,amscd}
\input amssym.def
\input amssym.tex

\newtheorem{thm}{Theorem}[section]

\newtheorem{cor}[thm]{Corollary}
\newtheorem{lem}[thm]{Lemma}
\newtheorem{prop}[thm]{Proposition}

\theoremstyle{remark}
\newtheorem{rek}[thm]{Remark}

\newcommand{\R}{\ensuremath{\mathbb{R}}}

\newcommand{\Z}{\ensuremath{\mathbb{Z}}}

\newcommand{\N}{\mathbb{N}}

\newcommand{\ud}{\mathrm{d}}
\newcommand{\dx}{\mathrm{d}x}

\newcommand{\abs}[1]{\left|#1\right|}

\title{Bounded step functions and factorial ratio sequences}
\author{Jason P. Bell \and Jonathan W. Bober}
\address{(J. Bell) Department of Mathematics, Simon Fraser University, Burnaby, BC, V5A 1S6, Canada}
\email{jpb@sfu.ca}
\address{(J. Bober) Department of Mathematics, University of Michigan, Ann Arbor, MI, 48109, USA}
\email{bober@umich.edu}
\date{March 27, 2008}
\begin{document}

\begin{abstract}
We study certain step functions whose nonnegativity is related to the
integrality of sequences of ratios of factorial products. In particular,
we obtain a lower bound for the mean square of such step functions which
allows us to give a restriction on when such a factorial ratio sequence
can be integral. Additionally, we note that this work has applications
to the classification of cyclic quotient singularities.
\end{abstract}

\maketitle

\renewcommand{\a}{\mathbf{a}}
\renewcommand{\b}{\mathbf{b}}
\newcommand{\balpha}{\mathbf{\alpha}}
\newcommand{\bbeta}{\mathbf{\beta}}

\newcommand{\f}[5]{ { }_{#1}F_{#2}\left(\begin{array}{c}#3 \\ #4\end{array};#5\right)}
\newcommand{\lcm}{\mathrm{lcm}}
\newcommand{\floor}[1]{\left\lfloor#1\right\rfloor}
\newcommand{\intersect}{\cap}
\newcommand{\frc}[1]{\left\{#1\right\}}

\newcommand{\rightsquare}{\begin{flushright}$\square$\end{flushright}}
\newcommand{\uab}{u_n(\a,\b)}
\newcommand{\uabfrac}{\frac{(a_1n)!(a_2n)!\cdots(a_Kn)!}{(b_1n)!(b_2n)!\cdots(b_Ln)!}}
\newcommand{\fab}{f(x;\a,\b)}
\newcommand{\gab}{g(n;\a,\b)}
\newcommand{\fabx}[1]{f(#1;\a,\b)}
\newcommand{\fabsum}{\sum_{k=1}^K \floor {a_kx} - \sum_{l=1}^L \floor{b_lx}}
\newcommand{\fabaltsum}{\sum_{l=1}^L \frc{b_l x} - \sum_{k=1}^K \frc{a_k x}}

\newcommand{\printlineno}{\hspace*{-.5in}\texttt{SOURCE LINE \#\the\inputlineno}}
\renewcommand{\printlineno}{}
\printlineno

\newcommand{\cmin}{c_{\min}}

\section{Introduction}

Let $u_n$ be the factorial ratio sequence
\[
    u_n = \uab  = \uabfrac,
\]
which depends on integer parameters $a_1, \ldots, a_K, b_1, \ldots, b_L \in \{1, 2, 3, \ldots\}$.
We consider the classification of parameters $\a = \{a_1, \ldots, a_K\}$ and $\b =  \{b_1, \ldots, b_L\}$
such that $u_n$ is an integer for all $n > 0$. It turns out that this is equivalent to studying the question of when
the step function
\[
    \fab = \fabsum
\]
is always nonnegative. In many interesting cases it is natural to restrict to $\a$ and $\b$ with
\[
    \sum_{k=1}^K a_k = \sum_{l = 1}^L b_l,
\]
so that $\fab$ is a periodic function with period $1$, and so that $\uab$ does not grow too fast.
For example, if these sums are equal and $\uab$ is always an integer, than $\a$ and $\b$ can easily
be used to give explicit elementary bounds on the prime counting function $\pi(x)$ via a standard
method (see \cite[Section 5.1]{shoup-ntb}, for example).

Two natural parameters associated with any such step function are the difference $h(L,K) = L-K$, which we will call the \emph{height},
and the sum $l(L,K) = L+K$, which we will call the \emph{length}.

If we restrict to the case of height $1$, then the situation is well understood
(see \cite{bober-factorial-ratios}, \cite{R-V-factorials}, \cite{vasyunin-step-functions}). Specifically, the set of pairs of tuples $(\a, \b)$ satisfying
\begin{equation}
\sum_{k=1}^K a_k = \sum_{l=1}^{K+1} b_l,
\end{equation}
\begin{equation}
\gcd(a_1, a_2, \ldots, a_K, b_1, b_2, \ldots, b_{K+1}) = 1,
\end{equation}
and
\begin{equation}
a_k \ne b_l \textrm{ for all } k,l,
\end{equation}
and such that $\uab$ is an integer for all $n$ is exactly known. There are three easily described two parameter infinite families of
such sequences and $52$ ``sporadic'' sequences which are easily listed.

A curious feature of this classification is that there are exactly two parameter sets of length $9$ and
none larger. It would be nice to have a ``simple'' reason to explain the fact that a function of
the form
\[
    f(x) = \sum_{k=1}^5 \floor{a_k x} - \sum_{l=1}^6 \floor{b_l x}
\]
cannot be positive for all $x$ if $\sum a_k = \sum b_l$ and $a_k \ne b_l$ for all $k,l$. Although the
classification in \cite{bober-factorial-ratios} proves that this is true, it does not seem to shed any light on exactly \emph{why} this is true.

More generally, if we fix the height $h(L,K)$ we may ask whether there is a reason that the length $l(L,K)$
must be small (in terms of the height) in order for $f(x)$ to have a chance of being always nonnegative.
We obtain here a result of this form as follows.

\begin{thm}\label{nonnegativity-bound}
Fix $L - K = D$. If
\begin{equation}\label{f-sum}
    f(x) = \fabsum
\end{equation}
where
\[
    \sum_{k=1}^K a_k = \sum_{l=1}^L b_l
\]
and $a_k \ne b_l$, $a_k, b_l \in \{1, 2, 3, \ldots \}$ for all $k,l$, and if
\[
    f(x) \ge 0
\]
for all $x$, then
\[
    K+L \ll D^2 (\log D)^2.
\]
\end{thm}

Similarly, we answer a conjecture of A. Borisov \cite{borisov-quotient-singularities} about
the boundedness of such step functions.

\begin{thm}\label{boundedness-bound}
If
\[
    f(x) = \sum_{k=1}^K \floor{a_kx}- \sum_{l=1}^L \floor{b_lx},
\]
where
\[
    \sum_{k=1}^K a_k = \sum_{l=1}^L b_l
\]
and $a_k \ne b_l$, $a_k, b_l \in \{1, 2, 3 \ldots \}$ for all $k,l$, and if
\[
    \abs{f(x)} \le A
\]
for all $x$, then
\[
    K+L \ll A^2 (\log A)^2.
\]
\end{thm}

\begin{rek}
These theorems will be stated and proved in a more explicit form in Section \ref{main-theorem-section},
and some comments will be made about the achieved bounds in Section \ref{remarks-section}.
\end{rek}

Concretely, Theorem \ref{nonnegativity-bound} says that if $L+K$ is large in terms of $L-K$, then
$f(x)$ must take negative values. Equivalently, this gives restrictions on when a factorial ratio sequence
can be an integer sequence. Similarly, Theorem \ref{boundedness-bound} says that if $L+K$ is large,
then $f(x)$ takes on large values.

The interest in Theorems \ref{nonnegativity-bound} and \ref{boundedness-bound} comes from
different areas of mathematics. On the one hand, an understanding of the nonnegativity of such step functions is
equivalent to an understanding of the integrality of sequences of ratios of factorial products.
Additionally,
after a change of variables, these step functions are some of the simplest functions
that show up in the Nyman--Beurling real variable reformulation of the Riemann hypothesis,
so it is possible that their study may shed some light on the subject---in fact, Vasyunin
\cite{vasyunin-step-functions} studied these functions from this perspective.
F. Rodriguez-Villegas \cite{rodriguez-villegas-hypergeometric-CY-manifolds} also has studied
the integrality of such step functions and relations to families of Calabi-Yau manifolds.
Moreover, A. Borisov \cite{borisov-quotient-singularities} has shown that in
some ways the classification of integral factorial ratios is equivalent to the classification
of cyclic quotient singularities having certain properties.

In fact, Borisov recasts Theorem \ref{boundedness-bound} as the following statement about cyclic quotient singularities.

\begin{prop}
Suppose $a \ge 0$ is any real number. Then for all large enough $d \in \N$, for all
but finitely many $(x_1, x_2, \ldots, x_d) \in T^d$ that define a cyclic
quotient singularity with Shukarov minimal log-discrepancy at least $d/2 - a$,
for some pair of indices $1 \le i < j \le d$ we have $x_i + x_j = 1$.
\end{prop}
\begin{proof}
The appears as \cite[Conjecture 3]{borisov-quotient-singularities}, where Borisov shows that
it follows from Theorem \ref{boundedness-bound}.
\end{proof}

Our theorems will actually be proved by giving a lower bound on the mean square of
a step function of the form \eqref{f-sum}. To give our bound we will notice that the Fourier coefficients of such
a function have some nice arithmetic properties. In fact, the Dirichlet series whose coefficients
are the Fourier coefficients of $f(x)$ is just the product of the Riemann $\zeta$-function and a
Dirichlet polynomial, so, combining Parseval's Theorem with a theorem of Carlson,
we are instead able to study the appropriate Dirichlet series on a vertical line.

Our original version of this paper proved Theorems \ref{nonnegativity-bound} and \ref{boundedness-bound}
with bounds that were exponential in $K-L$ and $A$. We are indebted to the K. Soundararajan for
an argument (namely, the approximation of the $\zeta$-function by a truncated Euler
product) which leads to much improved lower bounds in Theorem \ref{explicit-theorem} (and hence
Theorems \ref{nonnegativity-bound} and \ref{boundedness-bound}) and simplifies our exposition. After completing our original work
we learned that E. Bombieri and J. Bourgain \cite{bombieri-bourgain} also obtained our Theorem \ref{boundedness-bound}
through different methods.

\subsection{Notation}
$\floor{x}$ denotes the \emph{floor} of $x$, which is the largest integer less than
or equal to $x$. $\frc{x} = x - \floor{x}$ is the fractional part of $x$.
Also $e(x) := \exp(2\pi ix) := e^{2\pi ix}$. $\gamma$ denotes Euler's constant
\[
    \gamma := \lim_{N \rightarrow \infty} \left(\sum_{n = 1}^N \frac{1}{n} - \log(n)\right) \approx 0.577215664901533.
\]

\subsection{Acknowledgements}
We thank J. C. Lagarias for suggesting the problem that led to this
paper and for helpful comments on earlier drafts. We also thank A. Borisov for some comments
and for providing a preprint of \cite{borisov-quotient-singularities}, and
E. Bombieri for providing a manuscript with an alternate proof of our theorems.
\printlineno

\section{The connection between factorial ratios and step functions}\label{natural-function-section}
The connection that the step functions which we study have with factorial ratio sequences
comes from the following theorem.
\begin{prop}[Landau \cite{landau-factorials}] Let $a_{i,k}, b_{j,k} \in \Z_{\ge 0}, 1 \le i \le l, 1 \le j, \le k,
1 \le k \le r$ and let
\[
    A_i(x_1, x_2, \ldots, x_r) = \sum_{k=1}^r a_{i,k}x_k
\]
and
\[
    B_j(x_1, x_2, \ldots, x_r) = \sum_{k=1}^r b_{j,k}x_k.
\]
(That is, $A_i$ and $B_j$ are linear forms in $r$ variables with nonnegative integral
coefficients.) Then the factorial ratio
\[
    \frac{\prod_{i=1}^l A_i(x_1, x_2, \ldots, x_r)!}{\prod_{j=1}^k B_j(x_1, x_2, \ldots, x_r)!}
\]
is an integer for all $(x_1, \ldots, x_r) \in \Z_{\ge 0}^r$ if and only if the function
\[
F(y_1, \ldots, y_r) = \sum_{i=1}^l \floor{A_i(y_1, \ldots, y_r)} - 
        \sum_{i=1}^k \floor{B_j(y_1, \ldots, y_r)}
\]
is nonnegative for all $(y_1, \ldots, y_r) \in [0,1]^r$.
\end{prop}

A special case of this theorem is the following.

\begin{cor}\label{factorial-ratio-cor}Let 
\[
    \uab = \uabfrac
\]
Then $u_n$ is an integer for all $n$ if and only if
\[
    \fab = \fabsum
\]
is nonnegative for all $x$.
\end{cor}

In \cite{vasyunin-step-functions} Vasyunin considered functions of this type taking on only the values
$0$ and $1$. A simple generalization of a proposition of Vasyunin says that a necessary condition
that such a function take values only in the range $0 \ldots D$ is that $L - K = D$.

\printlineno

\begin{lem}\label{k-minus-l-condition}
Suppose that $f(x)$ is a function of the form
\[
    f(x) = \sum_{k=1}^K \floor{a_k x} - \sum_{l=1}^L \floor{b_l x}
\]
with $a_k, b_l \in \Z$, and that $f(x)$ is bounded. Then
$\sum_{k=1}^K a_k = \sum_{l=1}^L b_l$ and, for any $n$, there
exists some $x$ such that $f(x) = -n$ if and only if there exists
some $x'$ such that $f(x') = L-K + n$. In particular, $f(x)$ is nonnegative
if and only if the maximum value of $f$ is $L-K > 0$.
\end{lem}
\begin{proof}
The first assertion is clear, for if $\sum a_k \ne \sum b_l$, then
$f(x)$ is unbounded. Now we know that $f(x)$ is periodic with period $1$. Now, for
any $z$ that is not an integer we have $\floor{z} + \floor{-z} = -1$, so for
any $z$ for which none of $a_i z$, $b_j z$ is an integer, we have
\[
    f(z) + f(-z) = L-K,
\]
from which the assertion follows.
\end{proof}

\printlineno

\section{Proof of the main theorem}\label{main-theorem-section}
We study functions of the form
\[
    f(x;\a,\b) = \sum_{k=1}^K \floor{a_k x} - \sum_{l=1}^L \floor{b_l x}
\]
with
\[
    a_k, b_l \in \N
\]
subject to the conditions
\begin{equation}\label{sum-condition}
    \sum_{k=1}^K a_k = \sum_{l=1}^L b_l \textrm{ and } a_k \ne b_l \textrm{ for all } k,l.
\end{equation}
The first condition ensures that $\fab$ is periodic with period $1$ and the second is 
a nontriviality condition. Thus $\fab$ is a periodic integer valued step function that
has possible jump discontinuities at points $n/a_k$ or $n/b_l$ for $n \in \Z$.
Notice that, using the relation $\frc{x} = x - \floor{x}$ and the condition \eqref{sum-condition}, $\fab$ can be rewritten in
terms of combinations of fractional parts as
\[
    \fab = \fabaltsum,
\]
which is the form that we will generally use. In this section, we will give the following explicit
lower bound on the mean square of such step functions.
\begin{thm}\label{explicit-theorem}
Let
\[
    f(x) = \fabsum
\]
where
\[
    \sum_{k=1}^K a_k = \sum_{l=1}^L b_l
\]
and $a_k \ne b_l$, $a_k, b_l \in \{1, 2, 3, \ldots \}$ for all $k,l$, then for all $M > 285$
\begin{multline}\label{explicit-theorem-equation}
    \int_0^1 f(x)^2\dx \ge 
        \frac{(K - L)^2}{4} + \frac{(K+L)}{2\pi^2} \left[\frac{e^{-\gamma}}{(\log M)}\left(1 - \frac{1}{2 (\log M)^2}\right)\right]^2 \\
        - \frac{(K+L)^2}{2\pi^2} \left[ \frac{2}{(M-1)^{1/2}}\left(\frac{e^{-\gamma}}{(\log M)}\left(1 - \frac{1}{2 (\log M)^2}\right)\right)^{-1} + \frac{1}{M-1}\right].
\end{multline}
\end{thm}
\begin{rek}
In this theorem $M$ is a parameter that we will be free to choose to optimize this lower bound. We obtain
easily stated asymptotic results below by choosing $M = (K+L)^4$ and give examples of optimized
bounds that can be obtained in Section \ref{remarks-section}.
(The condition that $M > 285$ in the above statement comes from our use of Rosser and Schoenfeld's \cite{MR0137689} explicit
formulation of Merten's Theorem (see Lemma \ref{mertens-bound}).)
\end{rek}
Granting this theorem for a moment, we may now easily prove Theorems \ref{nonnegativity-bound} and \ref{boundedness-bound}.

\begin{proof}[Proof of Theorems \ref{nonnegativity-bound} and \ref{boundedness-bound}]
Take $M$ to be $(K + L)^4$. Then for fixed $K-L$, equation \eqref{explicit-theorem-equation} gives
\[
    \max_{x \in \R}\abs{f(x)}^2 \ge \int_0^1 f(x)^2 \dx \gg \frac{(K+L)}{(\log (K+L))^2}.
\]
Now let $N = K+L$ and $A = \max_{x \in \R} \abs{f(x)}$. Then we have
\[
    \frac{N}{(\log N)^2} \ll A^2.
\]
Additionally, since this implies that $A \gg N^{1/3}$, we may write
\[
    N \ll N \frac{(\log A)^2}{(\log N)^2} \ll A^2 (\log A)^2.
\]
This gives Theorem \ref{boundedness-bound}. Theorem \ref{nonnegativity-bound} follows
from Lemma \ref{k-minus-l-condition}, which said that $f$ is nonnegative if and only
if the maximum value of $f$ is $L - K$.
\end{proof}

To prove Theorem \ref{explicit-theorem}, we begin with the Fourier expansion for our step functions.
\begin{lem}\label{f-fourier-expansion}
Suppose that $f(x) = \sum_{l=1}^L \frc{b_l x} - \sum_{k=1}^K \frc{a_k x}$. Then the Fourier
expansion of $f$ is
\[
    f(x) = \frac{L-K}{2} + \frac{1}{2\pi i}\sum_{\substack{n\in\Z\\n\ne 0}}\frac{1}{n}\left[\sum_{a_k | n} a_k - \sum_{b_l | n} b_l\right]e(nx).
\]
That is, $\hat f(0) = (L-K)/2$ and for $n \ne 0$
\[
\hat f(n) = \frac{1}{2\pi in}\left[\sum_{a_k | n} a_k - \sum_{b_l | n} b_l\right].
\]
\end{lem}
\printlineno
\begin{proof}
The Fourier expansion for the fractional part of $x$ is
\begin{equation}\label{frac-x-fourier-expansion}
    \frc{x} = \frac{1}{2} - \frac{1}{2\pi i}\sum_{\substack{n \in \Z \\ n \ne 0}}\frac{e(nx)}{n}.
\end{equation}
Using this we can write
\[
    \fabaltsum = \frac{L-K}{2} - \frac{1}{2\pi i}\sum_{\substack{n \in \Z \\ n \ne 0}}
            \frac{1}{n}\left[\sum_{l=1}^L\ e(nb_lx) - \sum_{k=1}^K e(na_kx)\right].
\]
Changing the order of summation we get
\[
    \fabsum = \frac{L-K}{2} + \frac{1}{2\pi i}\left[
        \sum_{k=1}^K \sum_{\substack{n\in\Z\\n\ne 0}}\frac{e(na_kx)}{n} - 
        \sum_{l=1}^L \sum_{\substack{n\in\Z\\n\ne 0}}\frac{e(nb_lx)}{n}  
        \right].
\]
\printlineno

Now extracting the coefficient of $e(mx)$ in the sum we see that for $m \ne 0$ we have
\[
\hat f(m) = \frac{1}{2\pi i}\left[
    \sum_{\substack{n,a_k \\ na_k = m}} \frac{1}{n} -
    \sum_{\substack{n,b_l \\ nb_l = m}} \frac{1}{n}\right].
\]
Now the result follows on replacing the $n$ in the sum with $n = m/a_k$.
\end{proof}
\printlineno
\begin{rek}
From now on we set $D = L - K$. It is convenient to subtract off the first
Fourier coefficient of $f(x)$ and to consider
\[
    \int_0^1 \abs{f(x) - \frac{D}{2}}^2 \dx.
\]
On considering the Fourier expansion, it is easy to see that
\begin{equation}\label{misc-equation-4}
    \int_0^1 \abs{f(x)}^2 \dx = \int_0^1 \abs{f(x) - \frac{D}{2}}^2 \dx + \frac{D^2}{4},
\end{equation}
so it is simple to transfer a lower bound on one to a lower bound on the other.
Notice also that $\abs{\hat f(n)} = \abs{\hat f(-n)}$, so it follows from Parseval's theorem that
\begin{equation}\label{integral-fhat-relation}
    \int_0^1 \abs{f(x) - D/2}^2\dx = 2\sum_{n=1}^\infty \abs{\hat f(n)}^2.
\end{equation}
\end{rek}
\printlineno

\begin{rek}
We now notice a M\"obius inversion-type formula for the Fourier coefficients of $f(x) = \fab$.
Let
\begin{equation}\label{g-definition}
g(n) = \gab = \#\{a_k : a_k = n\} - \#\{b_l : b_l = n\}.
\end{equation}
Then from Lemma \ref{f-fourier-expansion} we see that for $n \ge 1$ we have
\[
    \hat f(n) = \frac{1}{2 \pi i} \sum_{d | n} \frac{dg(d)}{n}.
\]
Or, forming the Dirichlet series
\[
    G(s) = D(g,s) = \sum_{n=1}^\infty \frac{g(n)}{n^s}
\]
(note that $G$ is actually given by a finite sum)
and
\[
    F(s) = D(\hat f, s) = \sum_{n=1}^\infty \frac{\hat f(n)}{n^s},
\]
we have the relation
\begin{equation}\label{G-F-Relation}
    G(s)\zeta(s+1) = 2\pi iF(s),
\end{equation}
where $\zeta(s) = D(1,s) = \sum_{n=1}^\infty n^{-s}$ is the Riemann $\zeta$-function.
\end{rek}
\printlineno

To estimate the mean square of $f$ we use the following theorem of Carlson \cite{carlson}
to relate $\sum_{n=1}^\infty \abs{\hat f(n)}^2$ to a mean value of 
$\abs{G(it)}^2$.

\begin{prop}\label{carlson-theorem}
Let $f(s) = \sum_{n=1}^\infty a_n n^{-s}$ in some half-plane. Then if $f(s)$
is analytic and of finite order for $\sigma \ge \alpha$, and
\[
    \frac{1}{2T}\int_{-T}^T \abs{f(\alpha + it)}^2 \ud t
\]
is bounded as $T \rightarrow \infty$, then
\[
\lim_{T \rightarrow \infty} \frac{1}{2T} \int_{-T}^T \abs{f(\sigma + it)}^2 \ud t = \sum_{n=1}^\infty \frac{\abs{a_n}^2}{n^{2\sigma}}.
\]
\end{prop}
\begin{proof} See \cite[Section 9.51]{theory-of-functions}.\end{proof}

It follows immediately from this proposition and the relation \eqref{G-F-Relation}
that for all $\sigma > 0$ we have
\begin{equation}\label{misc-equation-A}
    \sum_{n = 1}^\infty \frac{\abs{\hat f(n)}}{n^{2\sigma}} = \lim_{T \rightarrow \infty} \int_{-T}^T \abs{G(\sigma + it)\zeta(1 + \sigma + it)}^2\ud t.
\end{equation}
Eventually we will let $\sigma$ tend to zero, but first we estimate the right hand side
of \eqref{misc-equation-A} using a truncated Euler product for the $\zeta$-function. Define
the main term
\[
    \zeta_M(s) = \prod_{p \le M} \left(1 - \frac{1}{p^s}\right)^{-1},
\]
and let $\zeta_R(s)$ be defined by $\zeta(s) = \zeta_M(s) + \zeta_R(s)$. We will
use the following lemmas which give bounds on the size of $\zeta_M$.

\begin{lem}\label{zeta-m-bound} If $\mathrm{Re}(s) \ge 1$ then
\begin{equation}
    \frac{1}{\zeta_M(1)} \le \abs{\zeta_M(s)} \le \zeta_M(1).
\end{equation}
\end{lem}
\begin{proof}
This is just a restatement of the fact that for all $s$ with real part
greater than $1$, and all integers $n > 0$, we have
\[
    \left(1 - \frac{1}{n}\right) \le \abs{1 - \frac{1}{n^s}}^{-1} \le \left(1 - \frac{1}{n}\right)^{-1}.
\]
\end{proof}
\begin{lem}[Effective Mertens' bound]\label{mertens-bound}
For all $M > 285$,
\begin{equation}
    \zeta_M(1) \le \frac{\log M}{e^{-\gamma}}\left(1 - \frac{1}{2(\log{M})^2}\right)^{-1}
\end{equation}
\end{lem}
\begin{proof}
See Rosser and Schoenfeld \cite[Theorem 7]{MR0137689}.
\end{proof}

\begin{proof}[Proof of Theorem \ref{explicit-theorem}]
We now write
\begin{multline*}
    \frac{1}{2T}\int_{-T}^T \abs{G(\sigma + it)\zeta(1 + \sigma + it)}^2 \ud t\\
     = \frac{1}{2T}\int_{-T}^T \abs{G(\sigma + it)\zeta_M(1 + \sigma + it)}^2 \ud t + E(M,T),
\end{multline*}
where
\begin{multline*}
E(M,T) = \frac{1}{2T}\int_{-T}^T\abs{G(\sigma + it)}\abs{\zeta_M(1 + \sigma + it) + \zeta_R(1 + \sigma + it)}^2 \ud t \\
       - \frac{1}{2T}\int_{-T}^T \abs{G(\sigma + it)\zeta_M(1 + \sigma + it)}^2 \ud t.
\end{multline*}
Applying the triangle inequality in the form $\abs{ \abs{A} - \abs{B}} \le \abs{A - B}$ to the above, we obtain
\begin{multline*}
\abs{E(M,T)} \le \frac{1}{2T}\int_{-T}^T\abs{G(\sigma + it)}^2\abs{\zeta_R(1 + \sigma + it)}^2 \ud t \\
     + \frac{1}{2T}\int_{-T}^T 2\abs{G(\sigma + it)^2\zeta_M(1 + \sigma + it)\zeta_R(1 + \sigma + it)} \ud t.
\end{multline*}
On inserting the bound for $\zeta_M$ from Lemma \ref{zeta-m-bound}, the bound $\abs{G(\sigma + it)} \le K + L$, and 
applying Cauchy's inequality to the integral $\int_{-T}^T \abs{\zeta_R(1 + \sigma + it)} \ud t$, we have
\begin{multline*}
\abs{E(M,T)} \le 2(K + L)^2 \zeta_M(1) \left[ \frac{1}{2T} \int_{-T}^T \abs{\zeta_R(1 + \sigma + it)}^2 \ud t\right]^{1/2} \\
     + (K + L)^2 \frac{1}{2T} \int_{-T}^T \abs{\zeta_R(1 + \sigma + it)}^2 \ud t.
\end{multline*}
Define $E(M) := \lim_{T \rightarrow \infty} E(M,T)$. Then
\begin{eqnarray*}
\abs{E(M)} &\le& 2(K + L)^2 \zeta_M(1) \left[\sum_{n \ge M} \frac{1}{n^{2 + 2\sigma}}\right]^{1/2} + (K + L)^2 \sum_{n \ge M}\frac{1}{n^{2 + 2\sigma}} \\
     &\le& 2(K + L)^2 \zeta_M(1) \left[\sum_{n \ge M} \frac{1}{n^{2}}\right]^{1/2} + (K + L)^2 \sum_{n \ge M}\frac{1}{n^{2}}.
\end{eqnarray*}
Inserting the bound
\begin{equation*}
\sum_{n \ge M}\frac{1}{n^{2}} \le \frac{1}{M - 1},
\end{equation*}
we get
\begin{equation*}
\abs{E(M)} \le \frac{2(K + L)^2 \zeta_M(1)}{(M - 1)^{1/2}}+ \frac{(K + L)^2}{M - 1}.
\end{equation*}

Now,
\begin{equation*}
    \frac{1}{2T}\int_{-T}^T \abs{F(\sigma + it)}^2 \ud t = \frac{1}{(2\pi)^2}\frac{1}{2T}\int_{-T}^T \abs{G(\sigma + it)\zeta_M(1 + \sigma + it)}^2 \ud t + \frac{E(M,T)}{4\pi^2},
\end{equation*}
so inserting the bound from our Lemma \ref{zeta-m-bound} and taking the limit as $T \rightarrow \infty$ we get
\begin{eqnarray*}
    \sum_{n = 1}^\infty \abs{\hat f(n)}^2 &\ge& \frac{1}{\zeta_M(1)^2}\frac{1}{(2\pi)^2} \sum_{n=1}^\infty \abs{g(n)}^2 + \frac{E(M)}{4\pi^2} \\
    &\ge& \frac{1}{\zeta_M(1)^2}\frac{1}{(2\pi)^2} \sum_{n=1}^\infty \abs{g(n)}^2 - \frac{\abs{E(M)}}{4\pi^2}.
\end{eqnarray*}
Putting in the bounds for $E(M)$, we get
\begin{equation*}
    \sum_{n = 1}^\infty \abs{\hat f(n)}^2 \ge \frac{1}{\zeta_M(1)^2}\frac{1}{(2\pi)^2} (K+L)
        - \frac{(K+L)^2}{4\pi^2}\left[\frac{2 \zeta_M(1)}{(M - 1)^{1/2}} + \frac{1}{M-1}\right].
\end{equation*}
And finally, putting in the bounds on $\zeta_M(1)$, we get
\begin{multline}\label{misc-equation-C}
    \sum_{n = 1}^\infty \abs{\hat f(n)}^2 \ge 
        \frac{(K+L)}{4\pi^2} \left[\frac{e^{-\gamma}}{(\log M)}\left(1 - \frac{1}{2 (\log M)^2}\right)\right]^2 \\
        - \frac{(K+L)^2}{4\pi^2} \left[ \frac{2}{(M-1)^{1/2}}\left(\frac{e^{-\gamma}}{(\log M)}\left(1 - \frac{1}{2 (\log M)^2}\right)\right)^{-1} + \frac{1}{M-1}\right].
\end{multline}
From this we obtain the following explicit lower bound on the mean square of $f$.
The theorem now follows immediately from \eqref{misc-equation-C} and fact that
\[\int_0^1 f(x)^2 \dx = 2\sum_{n=1}^\infty \abs{\hat f(n)}^2 + \frac{(K-L)^2}{4}.\]
\end{proof}
%
%
\section{Explicit bounds and remarks}\label{remarks-section}
\newcommand{\B}{\mathbf{B}}
Let $\B(D)$ denote the maximal number of terms that a step function of the form
\[
    \fab = \fabsum
\]
may have if $L - K = D$, $a_i, b_j \in \N$, $\sum a_i = \sum b_j$ and $a_i \ne b_j$ for
all $i, j$, and $f(x) \ge 0$ for all $x$. With this notation, Theorem \ref{nonnegativity-bound} says that
$\B(D)$ is finite, and that, moreover $\B(D) \ll D^2 \log D$. Using equation \ref{misc-equation-C}
we may compute some explicit bounds on $\B(D)$.

For example, when $D = 1$ such a function $f(x)$ that only takes on the values
$0$ and $1$ must have $\sum_{n=1}^\infty \abs{\hat f(n)}^2 = \frac{1}{4}$. 
If $K+L = 112371$ and $M = 112371^{4.96215}$, the right side of \ref{misc-equation-C}
is $\approx 0.250000802$.
Additionally, the right side of \eqref{misc-equation-C} is an increasing function of $K+L$ when
$M$ is set to be a fixed power of $K+L$, so it follows that $\B(1) < 112371$.

Similarly we obtain $\B(2) < 502827$ by taking $L+K = 502827$ and $M=502827^{4.6602}$, in which
case the right side of \eqref{misc-equation-C} is $\approx 1.00000138$, but it must be the case that
$\sum_{n=1}^\infty \abs{\hat f(b)}^2 \le 1$ if $f(x)$ takes values in precisely in the range
$\{0,1,2\}$.

These bounds are far from best possible.
In fact it is known that $\B(1) = 9$ (see \cite{bober-factorial-ratios}, which implies that $\B(D) \ge 9D$ for all $D$.
Computations by the second author suggest that it might be the case that $\B(2) = 18$.
Additionally, for the case of $D = 3$, a direct computer search by the second author has
checked that there are no functions with $29$ terms that take values only in the
range $\{0, 1, 2, 3\}$ when $\sum a_i = \sum b_j \le 60$. (However, the search
space would have to be greatly enlarged for this to be truly convincing evidence that $\B(3) = 27$.)
\printlineno

\bibliographystyle{amsplain}
\bibliography{/home/bober/math/electronic_papers/bib}

\end{document}